\documentclass[12pt]{amsart} 
\usepackage{color}
\usepackage{amssymb}
\usepackage{times}
\usepackage{amsmath}
\usepackage[pdftex]{graphicx}
\usepackage{hyperref}
\newcounter{sideremark}

\newtheorem{Thm}{Theorem}[section]
\newtheorem{Lem}[Thm]{Lemma}

\newtheorem{Def}[Thm]{Definition} 

\newenvironment{Pf}[1]
{\trivlist\item[]{\it #1\@. }}{\hspace*{\fill}$\Box$\endtrivlist}

\renewcommand{\marginpar}[1]{}
\catcode`\@=12

\def\Empty{}
\newcommand\oplabel[1]{
  \def\OpArg{#1} \ifx \OpArg\Empty {} \else
  	\label{#1}
  \fi}
		
\long\def\realfig#1#2#3#4{
\begin{figure}[htbp]
  \caption[#1]{#3}
  \centering
  \includegraphics[width=#4]{#2}
  \label{#1}
\end{figure}}

\newcommand{\comm}[1]{}

\renewcommand{\epsilon}{\varepsilon}

\renewcommand{\rho}{\varrho}

\begin{document}

\title{Wide short geodesic loops on closed Riemannian manifolds.}
\author{Regina Rotman}
\date{May 27, 2019}
\maketitle

\begin{abstract}
It is not known whether or not the lenth of the shortest periodic geodesic on 
a closed Riemannian manifold $M^n$ can be majorized by $c(n) vol^{ 1 \over n}$, or $\tilde{c}(n)d$, where $n$ is the dimension of $M^n$, $vol$ denotes the volume of $M^n$, and $d$ denotes its diameter. In this paper we will prove that for each $\epsilon >0$ one can find
such estimates for the length of a geodesic loop with with angle between 
$\pi-\epsilon$ and $\pi$ with an explicit constant that depends both on 
$n$ and $\epsilon$.


That is, let $\epsilon > 0$, and let 
$a = \lceil{ {1 \over {\sin ({\epsilon \over 2})}}} \rceil+1 $.  
We will prove that 
there exists a ``wide'' (i.e. with an angle that is wider than
$\pi-\epsilon$) geodesic loop
on $M^n$ of length at most $2n!a^nd$. 
We will also show that there exists a ``wide'' geodesic loop of 
length at most $2(n+1)!^2a^{(n+1)^3} FillRad \leq 2 \cdot n(n+1)!^2a^{(n+1)^3}
vol^{1 \over n}$. Here $FillRad$ is the Filling Radius of $M^n$. 


\end{abstract}


\section*{Introduction}


The two results in this paper are motivated by the following 35 year old question formulated  by M. Gromov in [G].  
Let $M^n$ be a closed
Riemannian manifold of dimension $n$ and volume $vol$. 
Gromov asked whether there exists a constant $c(n)$ such 
that the length of a shortest closed geodesic, $l(M^n)$, on $M^n$ 
is bounded above by 
$c(n) vol^{\frac{1}{n}}$. Similarly, one can ask if
there exists a constant $\tilde{c}(n)$, such 
that $l(M^n) \leq \tilde{c}(n)d$, where $d$ denotes the diameter of $M^n$. The answers to the above questions
are known for all closed Riemannian surfaces,
(see [BZ], [CK] for surveys of these results), including the most difficult case of a Riemannian $2$-sphere, for which the first upper bounds were established by C. B. Croke, (see [C]). Volume upper bounds were also proved by Gromov
in the special
case of essential manifolds, (see [G]). Finally, note that it  
easy to see that $l(M^n) \leq 2d$, when $M^n$ is not simply connected. 
This diameter upper bound does not have to hold in general, as it was shown by 
F. Balacheff, C. Croke and M. Katz who found an 
example of a Zoll $2$-sphere for which the length of a shortest periodic is greater than twice the diameter, (see [BCK]).

These are the only  known to us  curvature-free upper bounds for the length of a shortest periodic geodesic on closed Riemannian manifolds of dimension greater than two. The only   estimates for $l(M^n)$ in higher dimensions that involve curvature are (1) our old joint results with A. Nabutovsky,
[NR0] where we obtained the upper bound for $l(M^n)$  for the class of manifolds with the sectional curvature $K \geq k$, volume
$vol \geq v >0$ and diameter $d \leq D$,  and for the class of manifolds with 
$K \leq k$, and $vol \leq V$,  and (2) a recent theorem 
of N. Wu and Z. Zhu, (see [WuZ]), in which they established the existence of an upper bound for $l(M^n)$ for the class of $4$-manifolds with $vol \geq v$, $d \leq D$ and the Ricci curvature, $|Ric| \leq 3$. The result of Wu and Zhu uses recent theory of J. Cheger and A. Naber, (see [ChN]) of manifolds with two-sided bound on Ricci curvature. 


There is, however,  a number of curvature-free estimates for the length
of geodesic loops, starting with the result of S. Sabourau, [S] (in 
which he proved the first curvature-free upper bounds for the length of a shortest geodesic loop on an arbitrary Riemannian manifold) as well as for the length of stationary $1$-cycles, geodesic nets, and geodesic loops at each point of a manifold, (see [B], [NR1], [NR2], [R1], [R2], [R3], [R4]).
In fact, it is plausible  that  the curvature-free  upper bounds 
for the length of a shortest periodic geodesic do not exist in higher 
dimensions.
This would make the above results, as well as the result discussed in 
this paper optimal in some sense.

Let us begin with the following definitions, (see Definitions 0.1, 0.4 in [R4]):

\begin{Def} \label{DefA}

(a) A  minimal (or stationary) geodesic net $\Gamma$ 
is a multi-graph immersed into a Riemannian manifold $M^n$
satisfying the following two conditions: 

\noindent (1) Each edge of $\Gamma$ is a geodesic segment. Multiple edges between two vertices, integer multiples of the same edge, and edges between a point and itself, i.e. loops are allowed. 

\noindent (2) the sum of unit vectors at each vertex tangent to the edges
and directed away from this vertex is equal to zero. These vectors are counted with multiplicities of the respective edges. 

(b) If a minimal geodesic net $\Gamma$ has one vertex,
we will call it a (minimal, or stationary) geodesic flower, or a geodesic $m$-flower where $m$ is  the 
number of loops that comprise it counted with the multiplicities, (see fig. 
~\ref{flower}). Individual geodesic loops of the minimal geodesic flower will be called its petals. 
A minimal geodesic net $\Gamma$ that has at most $2$ vertices joined by
at most $m$ segments, (counted with multiplicities), will
be called  a minimal geodesic $m$-cage, (or just a geodesic
cage), (see fig. ~\ref{cages}). Thus, a geodesic $m$-flower can also be viewed as a geodesic $m$-cage.
\end{Def}

\noindent {\bf Remark.} Any immersion of a graph in $M^n$ will be 
called a net.  Nets with a single vertex will be called flowers, 
while nets with at most two vertices will be called cages. 

By Definition ~\ref{DefA} periodic  geodesics are minimal geodesic flowers, while singular geodesic loops are not. The stationarity condition implies that stationary geodesic nets, and, in particular, geodesic flowers are critical points of the (possibly weighted) length functional on the space of (multi-)graphs, (see [NR2]). 
That is a flower $F$ is a geodesic flower if and only if for any smooth vector field $v$ on a Riemannian manifold $M$, the functional $L(t)=$ {\it length}$(\Phi_t^v(F))$ has a critical point at $t=0$, where $\Phi_t^v$ is the $1$-parameter flow of diffeomorphisms of $v$. 

In this paper we will prove the following theorems:

\begin{Thm} \label{TheoremA}
Let $M^n$ be a closed Riemannian manifold of dimension $n$ and of
diameter $d$.  Let $q$ be the smallest integer such that 
$\pi_q(M^n) \neq \{ 0 \}$.
Then for any sufficiently small $\epsilon >0$ there exists a geodesic loop $\alpha$  with an  angle $\theta$, 
such that $\pi-\theta < \epsilon$, and such that its 
length $l(\alpha) \leq 2q!a^qd \leq 2n!a^nd$,
where $a=\lceil{{1 \over {\sin {\epsilon \over 2}}}} \rceil+1$. Moreover, this geodesic loop will be petal of a stationary geodesic flower. 
 
\end{Thm}

\begin{Thm} \label{TheoremB}
Let $M^n$ be a closed Riemannian manifold of dimension $n$ and of volume $vol$.
Then for every $\epsilon > 0$ there exists a geodesic loop $\alpha$ with an angle
$\theta$, such that $\pi-\theta <\epsilon$, and the length $l(\alpha) \leq
2 \cdot n (n+1)!^2 a^{(n+1)^3}vol^{1 \over n}$,
where $a=\lceil{{1 \over {\sin {\epsilon \over 2}}}} \rceil+1$ .  Moreover, this geodesic loop will be a
petal of  a stationary geodesic flower.

\end{Thm}




\realfig{cages}{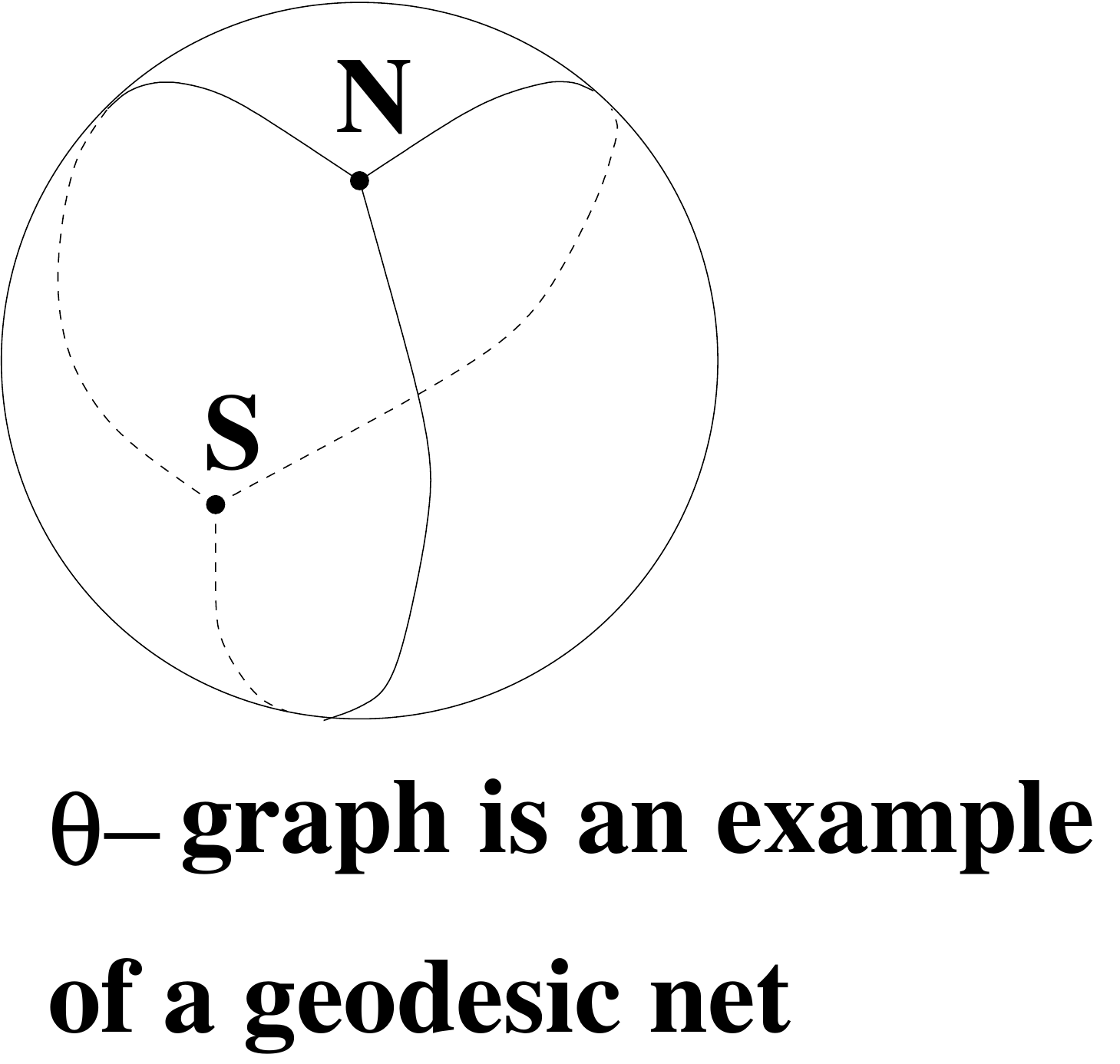}{A non-degenerate $3$-cage}{0.7\hsize}

\realfig{flower}{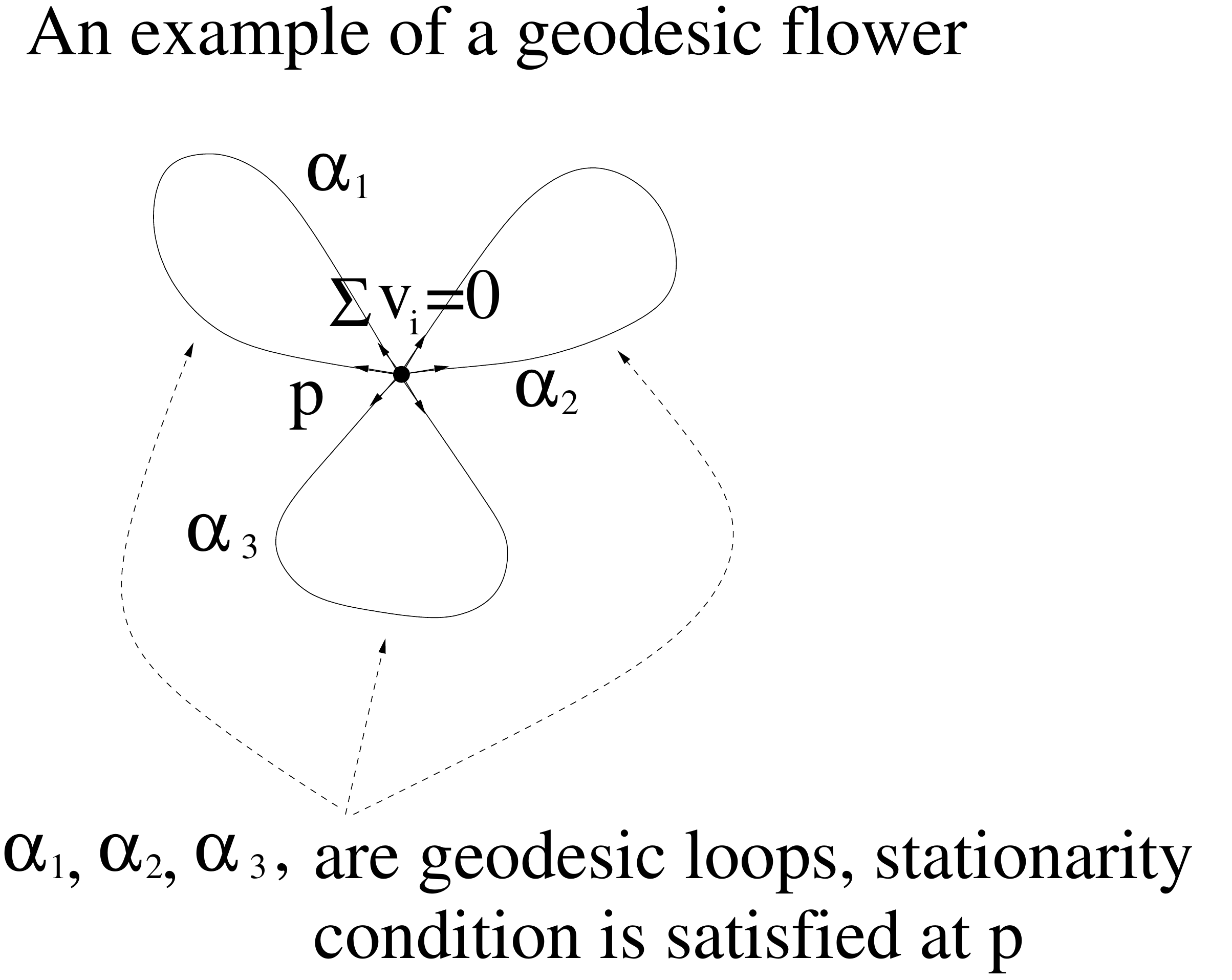}{A stationary geodesic flower}{0.7\hsize}
The first  idea behind the proofs of Theorems ~\ref{TheoremA}
and ~\ref{TheoremB} is a version of pseudo-filling argument pioneered by Gromov in 
[G], where one obtains a minimal object as an obstruction to filling a homologically non-trivial cycle in $M$. In our case, the wide geodesic loops will be obtained as a subset of a stationary geodesic flower the existence of which will be proven using such a technique. This filling technique is used in conjunction with two additional observations.

\noindent  (1) It is possible for the net to ``degenerate'' during the length shortening flow described in [NR0], whereby some of the edges can become smaller and eventually disappear, and the vertices can merge together. In view of this, one can introduce a weighted length shortening flow that will force 
all of the vertices in the net to merge together.

\noindent (2) The vertex stationarity condition implies that endowing one of the loops in 
the net with a large  weight will either force this loop to disappear or to become wide.


For example, the main idea of [R4] was  to define a weighted length functional on 
the space of nets such that its gradient-like  flow
will ``force'' some edges to shrink to a point, and nets  to 
degenerate into geodesic flowers.  

\begin{Def} 

\noindent (1) Let $\Gamma$ be a  net with edges
$e_1,...,e_i,..., e_k$. Then $L(\Gamma)=\Sigma_{i=1}^k m_i length(e_i)$,
where $m_i \in {\bf Z_+}$ and $length(e_i)$ is the length of $e_i$,
will be called a weighted length functional with weights $m_1,...,m_k$.  
In particular, given an integer $a \geq 3$, 
$L_a$ will denote the length functional with the weights $1,a, ..., a^{k-1}$.

\noindent (2) A net $N$ is stationary (or minimal, or critical)  with respect to a weighted length
functional $L$ with weights $m_i, i=1,...,k$, if for any one-parametric
smooth flow of diffeomorphisms $\psi_t, t=0$ is a critical point
of $f(t)=L(\psi_t(N))$.  

\end{Def}

The first variation of the length formula implies that if $N$ is critical with respect 
to a weighted length functional $L$, then (2) in the above definition implies that each edge of $N$ is a geodesic, and
the 
stationarity condition is satisfied at each vertex of $N$. Note that a critical point $N$ of the weighted length functional $L_a$ is a critical point of the regular length functional, if each edge of $N$ is taken with a multiplicity $m_i$.   

\realfig{stationarity}{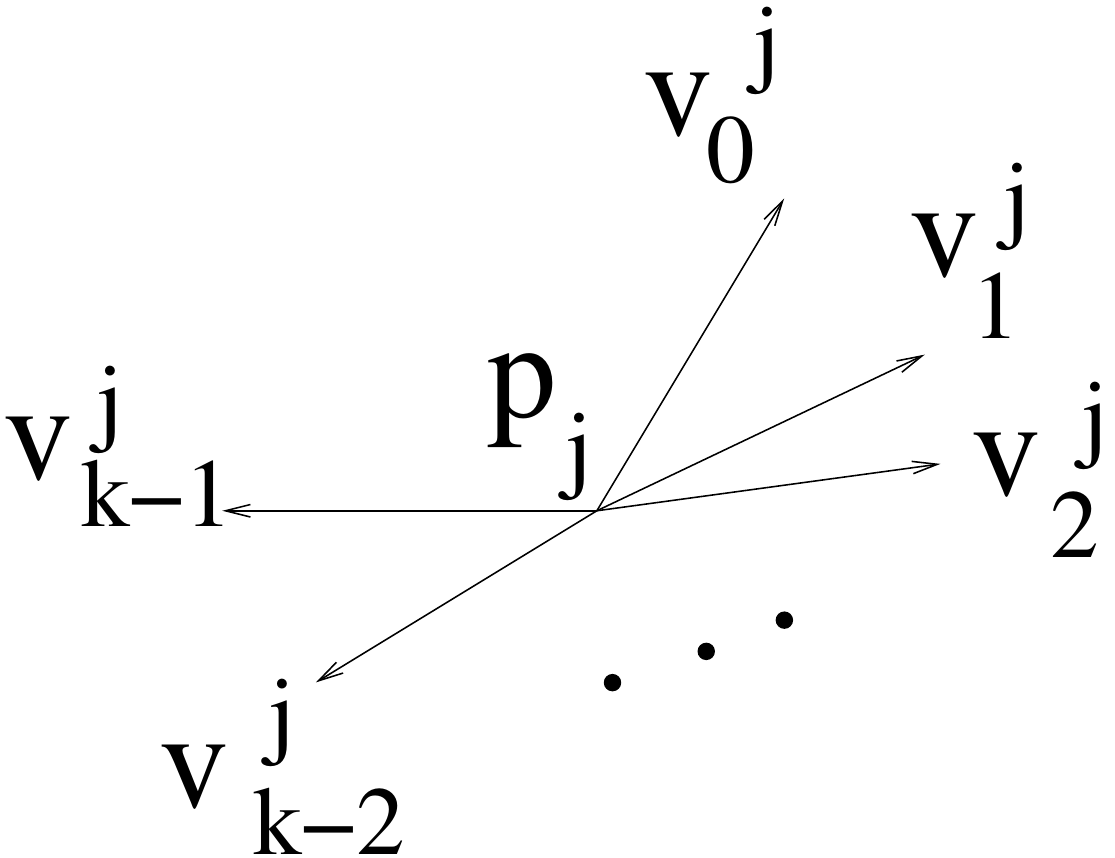}{Stationarity condition at $p_j$.}{0.7\hsize}


{\bf Example.} 
Let $\epsilon > 0$ be small enough so that  
$a=\lceil \frac{1}{\sin \frac{\epsilon}{2}} \rceil+1 \geq 3$.
Let $\Gamma$ be a $3$-cage with vertices $p_1$, $p_2$ and three edges
$e_1, e_2, e_3$. Note that it is possible that  $p_1 =p_2$ and  that the length of one or more edges of $\Gamma$ is zero. In the latter case, namely when 
one or more edges of $\Gamma$ is trivial $\Gamma$ will either be a closed curve, of a flower with two petals.    
Let $L_a(\Gamma)=length(e_1)+a \cdot length(e_2)+a^2 \cdot length(e_3)$ 
to be a weighted length functional applied to $\Gamma$ with the weights,
$1, a, a^2$. Suppose the weighted length shortening process converges to a non-trivial stationary geodesic cage.   

\realfig{theta}{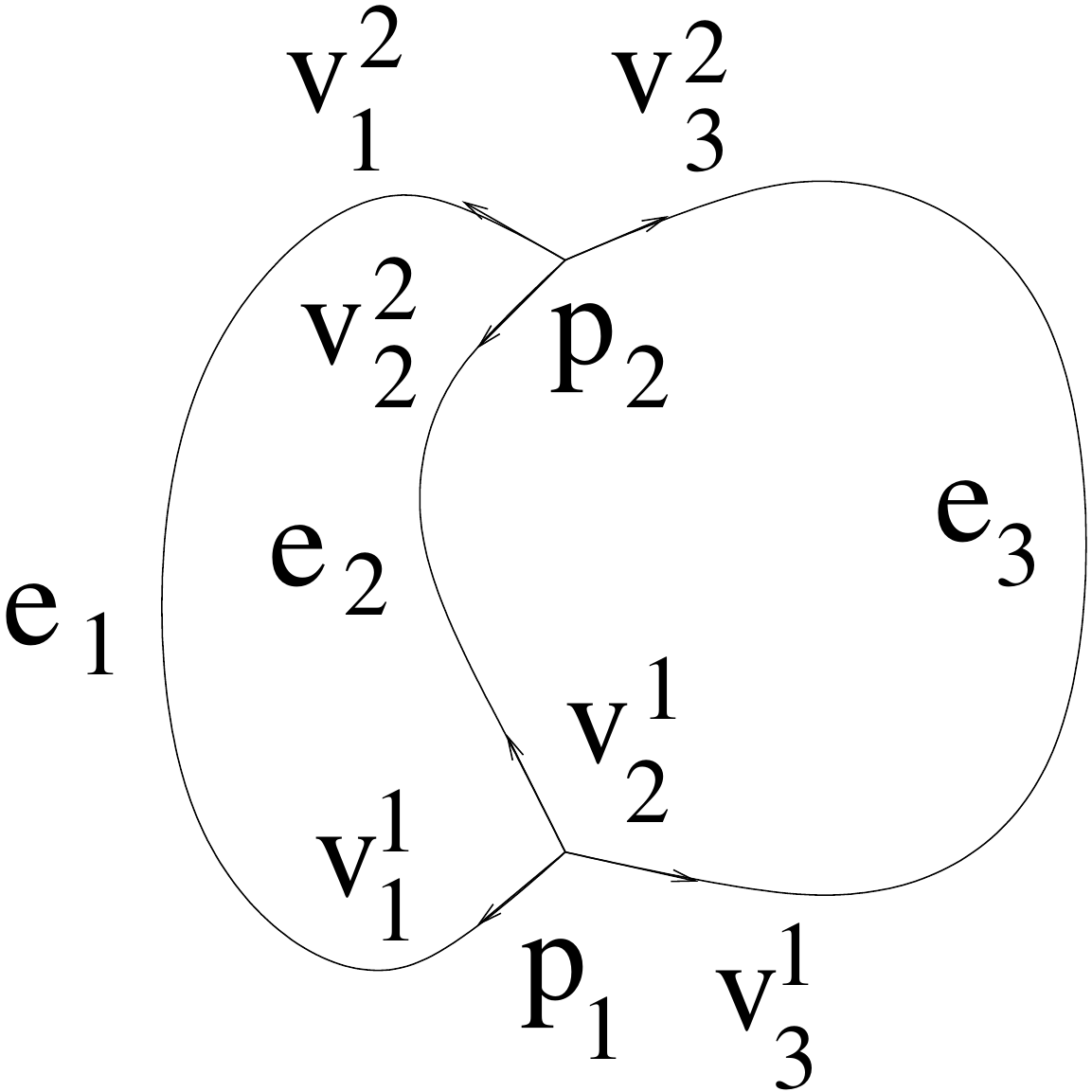}{Stationary $\theta$-graph.}{0.7\hsize}

Then a non-degenerate 
$3$-cage cannot be stationary. To show that suppose $\Gamma$ is not 
degenerate. Let $v_i^j$,  be the  unit vector tangent to $e_i$ at the point
$p_j, j\in \{1,2\}$, (see fig. ~\ref{theta}). Then the stationarity condition implies that
for both $j=1,2$, $v_1^j+av_2^j+a^2v_3^j=0$. These conditions obviously 
cannot be satisfied as $a^2 > 1+a$. Therefore, if $\Gamma$ is critical
then at least one of the $e_i$'s has to shrink to   a point and
the cage degenerates into a flower. 

We will next show that one of the ``petals'' of this flower has to 
be a geodesic loop with a wide angle, i.e. with an angle $\theta$, such 
that $\pi - \theta < \epsilon$. There are four possibilities to consider. 

{\bf Case 1: Two edges of $\Gamma$ shrink to a point.} In this case 
the stationary geodesic flower will be a periodic geodesic, which 
can be viewed as the geodesic loop with an angle $\pi$. The length of
the periodic geodesic will be at most $l_1+al_2+a^2l_3$, where $l_i$ is 
the length of $e_i, i \in \{1,2,3\}$.

\realfig{eight}{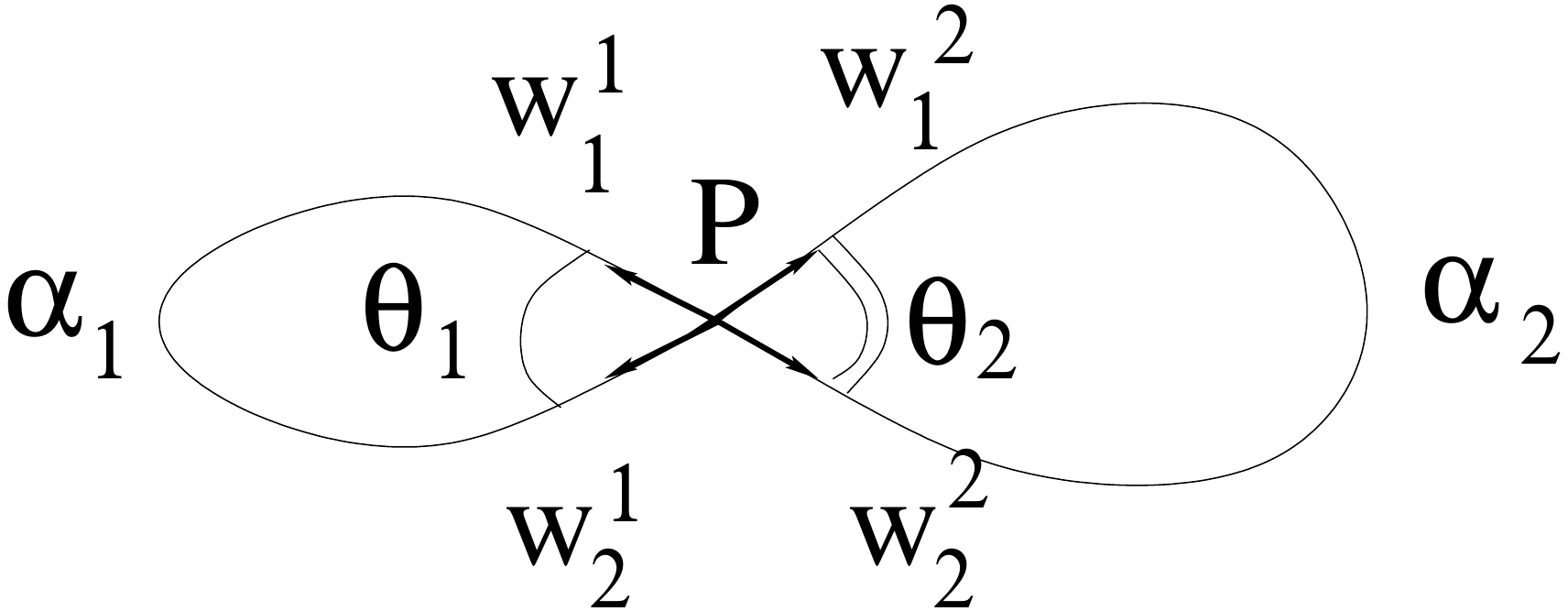}{Stationary figure $8$.}{0.7\hsize}

{\bf Case 2: $e_3$ shrinks to a point.} In this case the stationary
geodesic flower will be the ``figure 8'', (see fig. ~\ref{eight}). It will have two loops $\alpha_1$ obtained from the deformation of $e_1$ and $\alpha_2$ that is obtained from the deformation of $e_2$, and one
vertex $P$. Let $\theta_i$ be an angle corresponding to $\alpha_i$, for $i \in \{1,2\}$. The stationarity condition implies that the bisectors of $\theta_1$ and
$\theta_2$ are the same. Let $w_i^j$ be the unit vectors tangent to $\alpha_i$ at
$P$, $j \in \{1,2\}$. Then the lengths of the projections of $w_i^j$ onto the bisector will 
be $\cos \frac{\theta_i}{2}$. The stationarity condition then implies the following
equation $ \cos \frac{\theta_1}{2} = a \cos \frac{\theta_2}{2}$. Thus, 
$a \sin \frac{(\pi-\theta_2)}{2} \leq 1$, It follows that 
$\sin \frac{(\pi-\theta_2)}{2} \leq \frac{1}{a} < \sin \frac{\epsilon}{2}$.
Therefore, $\pi -\theta_2 < \epsilon$. 

{\bf Case 3: $e_2$ shrinks to a point.} 
Once again, the stationary geodesic net will be the ``figure 8'' with 
loops $\alpha_1, \alpha_3$ correpsonding to the edges $e_1, e_3$. with 
the vertex angles $\theta_1, \theta_3$. The bisectors of the angles will
be the same. Let $w_i^j$ be the unit vectors tangent to $\alpha_i, i \in \{1,3\}$, and $j \in \{1,2\}$. The lengths of the projections of the vectors onto the bisector
will be $\cos \theta_i$, and the stationarity condition will imply 
that $\cos \frac{\theta_1}{2}=a^2 \cos \frac{\theta_3}{2}=a^2 \sin \frac{(\pi-\theta_3)}{2}$. Therefore,  $\sin \frac{(\pi-\theta_3)}{2} \leq \frac{1}{a^2} < 
(\sin \frac{\epsilon}{2})^2 < \sin \frac{\epsilon}{2}$. Therefore, 
$\pi - \theta_3 < \epsilon$. 

{\bf Case 4: $e_1$ shrinks to a point.} The case is done in a similar
fashion.

We then combine
this idea with the techniques of [R4] that will be explained in the 
next section. Let us describe this more
formally.

To prove Theorem 0.3 we will use the weighted length functional $L_a$ with the weights $1, a, a^2, ..., a^{m-1}$, where $m$ is the number of edges in the $m$-cage. 
Recall that $a$ is a constant that depends on $\epsilon$, and that increases unboundedly as $\epsilon$ approaches $0$. Likewise, in order to prove 
Theorem 0.4, we will consider nets that correspond to the $1$-skeleton of
an $m$-simplex, and will use the weighted length functional with 
with weights $1, a, a^2, ..., a^k$, where $k={(m+1)(m+2) \over 2}$, i.e.
the number of edges in the $1$-skeleton.

The  non-trivial critical points that correspond to the above weighted length
functionals
are  bouquets of geodesic loops of total length 
$\leq \tilde{c}(n)d$ and of total length $\leq c(n) vol^{1 \over n}$ respectively, 
where $d$ is the diameter and $vol$ is the volume of $M^n$. Moreover, 
one of the loops in this net, namely the loop corresponding to the edge
with the highest weight,  among all of the edges of nonzero length, will have an angle that will approach $\pi$ as $a$ becomes large. We will prove that there always exists such a critical point, 
using techniques of [R4], modified Gromov's extension technique appearing in 
[G],  and
the idea illustrated by the above example.  

\section{The proof of Theorem ~\ref{TheoremA}.}

First we will need to prove the following lemma.

\begin{Lem} \label{Lemmadiameter} Let $\epsilon >0$ be given. 
Let $a = \max \{\lceil \frac{1}{\sin \frac{\epsilon}{2}} \rceil +1, 3\}$. 
Let $\Gamma$ be a geodesic cage with 
edges $e_i, i \in \{0,...,(k-1)\}$ each taken with the multiplicity
$a^{i}$.  Then $\Gamma$ is not stationary with respect to $L_a$, unless it is a geodesic flower. 
Moreover, 
let $F$ be a geodesic flower with the geodesic loops (petals) $e_i, i=0,...,k-1$, 
(some of them possibly
trivial) with multiplicities $a^i$. Then if $F$ is stationary, one 
of the petals has an angle  $\theta > \pi -\epsilon$.
\end{Lem}

\begin{Pf}{Proof} Suppose $\Gamma$ is stationary with respect to $L_a$. 
Let us denote the two 
vertices of the geodesic cage $\Gamma$ as $p_1, p_2$. Assume that $p_1 \neq p_2$. Let $v_i^j$ be 
the unit vector tangent to the edge $e_i$ at the vertex $p_j$, and diverging
from $p_j$, (see fig. ~\ref{stationarity}).  The stationarity condition implies that 
$\Sigma_{i=0}^{k-1}a^iv_i^j=0$, which can be restated as 
$a^{k-1}v_{k-1}^j=-\Sigma_{i=0}^{k-2}a^i v_i^j$, for $j \in \{1,2\}$. The length of the  
projection of $v_i^j$ onto the line passing through $v^j_{k-1}$ is  
$\cos \theta_i^j$, where $\theta_i^j$ is an angle that $v^j_i$ is making with
$-v^j_{k-1}$. Then the stationarity conditions becomes
$a^{k-1}=\Sigma_{i=0}^{k-2}a^i \cos \theta_i^j \leq \Sigma_{i=0}^{k-2}a^i=
\frac{a^{k-1}-1}{a-1}$, which is a contradiction. 

Next let us consider a geodesic flower $F$. Let $P$ be its vertex. 
We will denote the two unit vectors that are tangent to the non-constant loop 
$e_i$ and diverging away from $P$ as $w_i^1, w_i^2$. In the case when $e_i$ is a point, we will
let $w_i^j=0, j \in \{1, 2 \}$. Then the stationarity condition at $P$ implies that $\Sigma_{j=1}^2 \Sigma_{i=0}^{k-1} a^i w_i^j = 0$.  Let $s$ be the maximum index
for which the loop $e_j$ is not trivial.
Let $V = \Sigma_{j=1}^2\Sigma_i^{s-1} a^i w_i^j$. Then $a^s(w_s^1+w_s^2)=-V$.
The length of the  projection of $w_i^j$ onto the line through
$-V, i \in \{1,..., s-1\}$, where $\theta_i^j$ is the angle between $V$ and 
$w_i^j$ is $\cos \theta_i^j$. Also, let $\cos \frac{\theta}{2}$ be the projections of $w_s^j, j \in \{1, 2\}$
onto the line through $V$. Here $\frac{\theta}{2}$ is an angle that $w_s^j$ makes with
$-V$. Then the stationarity condition implies
$2a^s \cos \frac{\theta}{2} = \Sigma_{j=1}^2 \Sigma_{i=0}^{s-1} a^i \cos \theta_i^j \leq
\Sigma_{j=1}^2 \Sigma_{i=0}^{s-1}a^i \leq \frac{2(a^s-1)}{a-1}$. This implies 
that $\sin \frac{\pi-\theta}{2}=\cos \frac{\theta}{2} < \frac{1}{a-1} \leq 
\sin \frac{\epsilon}{2}$. Therefore, $\pi-\theta < \epsilon$.
\end{Pf}

\begin{Pf}{Proof of Theorem ~\ref{TheoremA}} Let $\epsilon >0$ be given such that 
$\lceil \frac{1}{\sin \frac{\epsilon}{2}} \rceil +1 \geq 3$.
The proof is by contradiction.
Let $f: S^q \longrightarrow M^n$ be a 
non-contractible map from the round sphere $S^q$ of dimension $q$ into $M^n$. 
Let us endow $S^q$ with a sufficiently fine triangulation, so that the the
maximal diameter of the simplices on $M^n$ in the triangulation induced by 
$f$  is at most some small positive $\delta$. Let us consider 
the disc $D^{q+1}$, such that $\partial D^{q+1} = S^q$. We will triangulate
it as a cone over $S^q$. 

Now consider the space $\Gamma M^n$ of the $m$-cages on $M^n$ together 
with the weighted 
length functional with weights: $1, a, a^2, ..., a^{m-1}$. Note that it is enough
to prove that there exists a critical point of ``small'' ( i. e. bounded in 
terms of the diameter $d$ of $M^n$) but non-zero length on $\Gamma M^n$. Indeed, by 
Lemma 1.1, this critical point will have to be a stationary geodesic flower, 
and one of the loops that comprise it will have an angle that is $\epsilon$-close to $\pi$. 

Thus, assume that all of the critical points that correspond to 
the above weighted length functional are longer than the bound specified in the conclusion of the theorem. We will then show that we can extend the map $f$ to 
$D^{q+1}$, which would contradict the fact that $f:S^q \longrightarrow M^n$ is 
not contractible. Now let us enumerate the edges in the $1$-skeleton of $D^{q+1}$. 
The pseudo-extension procedure will be done by induction on the  
skeleta of $D^{q+1}$.  Note that the map can be naturally extended to the 
 center of the disc by assigning to it
an arbitrary point $p$ in $M^n$. Likewise, the extension to the $1$-skeleton is
trivially accomplished by mapping
 the edges 
to minimal geodesic segments that connect this point with the corresponding
vertices of the triangulation of the image sphere.  The rest of the extension is done via an inductive bootstrap procedure, and amounts to ``filling'' $m$-cages by $m$-discs for all values of
$m \leq q+1$. A similar procedure was used in [R1] and [R4].
One proves the base of induction by contracting $2$-cages, i.e. closed
curves, by the usual Birkhoff curve shortening process, assuming there are no ``short'' periodic geodesis, i.e. geodesics of length at most $2d$. These homotopies generate $2$-disks that fill the $2$-cages.

Now assume that we have extended our map to the
$k$-skeleton,
we will now describe how to extend $f$  to the $(k+1)$-skeleton of $D^{q+1}$. In
order to do that we will show how to extend $f$ to each $(k+1)$-dimensional simplex
of $D^{q+1}$. It will be done by 
``filling'' $(k+1)$-cages by $(k+1)$-dimensional discs.
Suppose we want to extend $f$ to simplex $\sigma_{k+1}$. Consider
its boundary.
It consists of $k+2$ $k$-dimensional simplices.  One of the simplices in 
the boundary is a simplex of $S^q$. Recall that the diameter of its image 
under $f$ is at most some small number $\delta$, which implies that it 
can be effectively treated as a point, $q$ (see the 
remark on page 13 in [R1]. Since the simplex is contractable, one can 
continuously deform it to a point, while stretching each edge by at most 
$\delta$).
Consider an ordered  $(k+1)$-cage
$C$ that consists of the points $p, q$ and the geodesic edges that we can 
arrange into a sequence $e_0, e_1,..., e_k$ starting with an edge with the smallest index, and
ending with the edge with the largest simplex. Next we will apply the weighted
length shortening process, where the weight corresponding to $e_j$ will be
$a^j$. Let us now consider a weighted-length shortening process on 
$\Gamma M^n$. The construction of this flow is analogous to the Birkhoff 
Curve Shortening Flow. Likewise, it can be shown that in the absence of non-trivial 
critical points it is possible to deform $\Gamma M^n$ to the subspace that
consists of the constant cages, i. e. the cages, in which both vertices
coincide and all of the edges have zero lengths, (see [C] for
the detailed description of the Birkhoff curve shortening process). 
This flow was formally defined for $1$-cycles in a much more general setting 
in [NR1] and [R1]. The weighted length shortening flow is no different then 
the length shortening flow applied to the multi-graphs, where each edge is 
allowed to be taken with multiplicities.
The deformation of the edge during the length shortening is uniquely defined. 
Thus,  each copy of the multiple edge will be deformed in the same way, and the splitting of the edges is not possible.  
(See  Section 3  of [NR1] for the detailed description of the length shotening 
flow for nets).

If there are no critical geodesic flowers corresponding to the length functional $L_a$, then any cage $C$ that corresponds to the $1$-skeleton of $(k+1)$-simplex in $M^n$ obtained on the prior step of the induction can be 
contracted to a point that we will denote $x$ along a $1$-parameter
family of cages $C^{k+1}_\tau$, $\tau \in
[0,1]$ of smaller weighted length.  We can next construct
a $1$-parameter family of spheres $S^k_\tau$ of dimension $k$ corresponding to 
$C^{k+1}_\tau$, where  $S^k_0$ will be
the image of the boundary of the given simplex, and $S^k_1=\{x\}$.
This $1$-parameter family of spheres generates a $(k+1)$-dimensional disc.
Spheres are constructed by the procedure of ``filling'' cages 
$C^{k+1}_\tau$ at each $\tau$ first described in [R1].   For each $\tau \in [0,1]$, consider $(k+1)$ of $k$-''subcages'' of $C^{k+1}_\tau$. That is, we are considering $k$-cages obtained from $C^{k+1}_\tau$ by deleting from it one of the edges,
i.e. we are looking at all of the  $1$-parameter families of $k$-tuples $(e_0)_\tau,...,(\hat{e}_j)_\tau,...,(e_k)_\tau$
obtained from the original $(k+1)$-cage by removing from it one of the curves.  By induction assumption for all $\tau \in [0,1]$,
each of these subcages can be ``filled'' by discs of dimension $k$, and 
moreover, the resulting discs will change continuously with respect to $\tau$.

Next glue these $(k+1)$ $k$-dimensional discs as in the boundary of 
$(k+1)$-simplex to obtain $S^k_\tau$. Recall that the last $(k+2)$nd disk is 
a point. Moreover, 
this process is continuous with respect to $C^{k+1}_\tau$. 
This one-parameter family
of spheres generates the desired  $(k+1)$-dimensional disc that can
be used to extend $f$, (see fig. ~\ref{fig1}, which depicts how contracting 
$3$-cage generates a $1$-parameter family of $2$-spheres).  

We will next prove the length upper bound. 
Consider the $1$-parameter family $C^{k+1}_\tau$ obtained during the 
weighted length shortening flow of $C^{k+1}_0$ with the aformentioned weights
$1, a,...,a^k$. Let us denote the maximal total weighted length of $C^{k+1}_\tau$ over all $\tau \in [0,1]$  as 
$l_{k+1}$. Let us consider an ordered sequence of edges of 
$C^{k+1}_\tau$. 
Then the length of $(e_i)_\tau$ is at most ${l_{k+1} \over a^{i-1}}$  
for each $\tau \in [0,1]$ and for all $i \in \{0,...,k\}$. 
For each $\tau \in [0,1]$ let us consider all of the $k$-subcages of 
$C^{k+1}_\tau$. Their total length, of course, does not exceed the maximal length of the $(k+1)$- cage, i. e.  
$l_{k+1}$. Let us consider the $k$-subcage, denoted as $C^k$ of $C^{k+1}_\tau$ that is comprised
of $(e_0)_\tau,...,(e_{k-1})_\tau$. From the estimates of the individual segments, we can see that its total length can be  at most 
$l_{k+1}+{l_{k+1} \over a}+...+{l_{k+1} \over a^{k-1}}$. This is the 
upper bound valid for 
all of the other subcages of $C^{k+1}_\tau$ for all $\tau$. 
As we apply the weighted length shortening  with the coefficients $1, a, ..., a^{k-1}$ to $C^k$, we will obtain a 
$1$-parameter family $C^k_s, s \in [0,1]$, where $C^k_0=C^k$. Note that
the maximal length of the $k$-subcages resulting from this flow will be
at most
$l_{k+1}+a ({l_{k+1} \over a})+...
+a^{k-1}({l_{k+1} \over a^{k-1}}) =k \times l_{k+1}$. 
Let $l_k$ denote the maximal length over all of the (2-parameter) family of
the $k$-subcages. (One parameter is $\tau \in [0,1]$; the second parameter is 
the discrete parameter describing the choice of $k$ geodesics out of $(k+1)$ segments in the cage).  Then we obtain the following recurrent relation between
$l_{k+1}$ and $l_k$, $l_k \leq k \cdot l_{k+1}$.  This relation will hold for all 
$k=q, q-1,...,2$.
Now note that originally the length of each edge is bounded by the diameter
of the manifold $d$. Thus, $l_{q+1} \leq d(1+a+...a^q) ={d(a^{q+1}-1) \over a-1}
\leq 2da^q$. Putting this together we see that
$l_2 \leq 2q!a^qd$.

\end{Pf}
  
{\bf Example.} Let us separately present the proof of Theorem ~\ref{TheoremA} in a simple 
case of  $q=2$.
Let $f: S^2 \longrightarrow M^n$ be non-contractible. Let $a \geq 3$ be given.  The proof
will be by contradiction.  We will show that in the absence of stationary 
``figure 8'', (a flower with two geodesic loops) 
on $M^n$ of length $\leq 4a^2d$ that contains a loop with 
an angle that is $\epsilon$-close to $\pi$,  we 
can extend  $f$ to $D^3$ triangulated as a cone over $S^2$. Recall that
$S^2$  is 
triangulated in such a way that the image of a simplex under $f$ has
diameter at most $\delta$ for some small $\delta$ that eventually will
approach $0$.
The extension will constructed by induction on the skeleta of $S^2$.

{\bf 0-skeleton} of $D^3$ consists of the vertices of the triangulation of $S^2$ and of
the 
center of the disc $\tilde{p}$.  Thus, we extend to the $0$-skeleton by mapping
$\tilde{p}$ to  
an arbitrary point $p \in M^n$. 

To extend to the
{\bf $1$-skeleton}
of $D^3$ assign to an arbitrary $1$-edge of the form $[\tilde{p}, \tilde{v}_i]$, where 
$\tilde{v}_i$ is a vertex of the triangulation of $S^2$ a shortest geodesic segment $[p, v_i]$, where $v_i=f(\tilde{v}_i)$ of length at most $\leq d$ that connects $p$ with $v_i$.  

We extend to the 
{\bf $2$-skeleton} by assigning to an arbitrary $2$-simplex of
the form $[\tilde{p}, \tilde{v}_{i_1}, \tilde{v}_{i_2}]$
the disk generated by the contraction of the image of its boundary to some
point. The
boundary is mapped to a closed curve of length $\leq 2d+\delta$. Thus, in 
the absence of short periodic geodesics it is contractible via the Birkhoff curve shortening process.


Finally, let us extend to the {\bf $3$-skeleton}.  Consider an arbitrary $3$-simplex $[\tilde{p}, \tilde{v}_{i_1},\tilde{v}_{i_2},
\tilde{v}_{i_3}]$.  The image of its boundary was defined on the previous step of the induction. It is 
a $2$-sphere glued from four $2$-simplices. One of the four simplices comes from the triangulation of $S^2$, and can, therefore,   be made arbitrarily
small. One can, thus, contract it to $q$ over itself. This allows us to simply
treat it as a point $q$. So, one can view the sphere, i. e. the image of 
the boundary of $[\tilde{p}, \tilde{v}_{i_1}, \tilde{v}_{i_2}, \tilde{v}_{i_3}]$ as being formed by connecting two points $p$ and $q$ by
geodesic segments $e_1, e_2, e_3$ and then contracting each pair of closed 
curves to a point.  Note also, that this construction provides us
with a natural cell decomposition of this sphere into two
$0$-cells: $p, q$, three $1$-cells: $e_1, e_2, e_3$ and
three $2$-cells.  Consider the net that corresponds to the 
 $1$-skeleton of this sphere under this
decomposition.  
 Let us shorten  the weighted length $l(e_1) +
al(e_2)+a^2l(e_3)$. Under the above weighted length shortening, the net will either contract to a point, or will converge to a critical point, which 
can either be a stationary
figure $8$ or a periodic geodesic of length $\leq 4a^2d$.  Thus, assuming 
there are no critical points, the net can only converge to a point,
(see
fig. ~\ref{fig1}). This net shortening generates 
a $1$-parameter family of nets that we will denote
$C_\tau$. We can extend $C_\tau$ to one parameter family of 
$2$-spheres, $S^2_\tau$. $S^2_0$ will be the original sphere 
that corresponds to the image of the boundary of the $3$-simplex that we are trying to ``fill'', while $S^2_1$ will be  a point. We thus will obtain the
$3$-disk.
It now remains to show how to construct $S^2_\tau$.  For each $\tau \in [0,1]$ we consider
three pairs of curves and contract them to a point without the length increase.
At some point that the length of one or two 
segments will decrease to zero and segments themselves will degenerate to 
points.  
However we can still consider three pairs of curves, 
where one of the curves in two pairs will be a constant curve.  
We then fill each of these three pairs of curves by discs as we did
above when we were extending to the $2$-skeleton, i. e. using the curve
shortening process.  These $3$ $2$-discs glued together form 
$S^2_\tau$. 
Thus, if there is no geodesic stationary figure 8 (or periodic geodesic)
of length $\leq 4d$ that is a critical point with respect to $L_a$
then we can extend our map $f$ to the $3$-skeleton of $D^3$,  
reaching a contradiction.

\realfig{fig1}{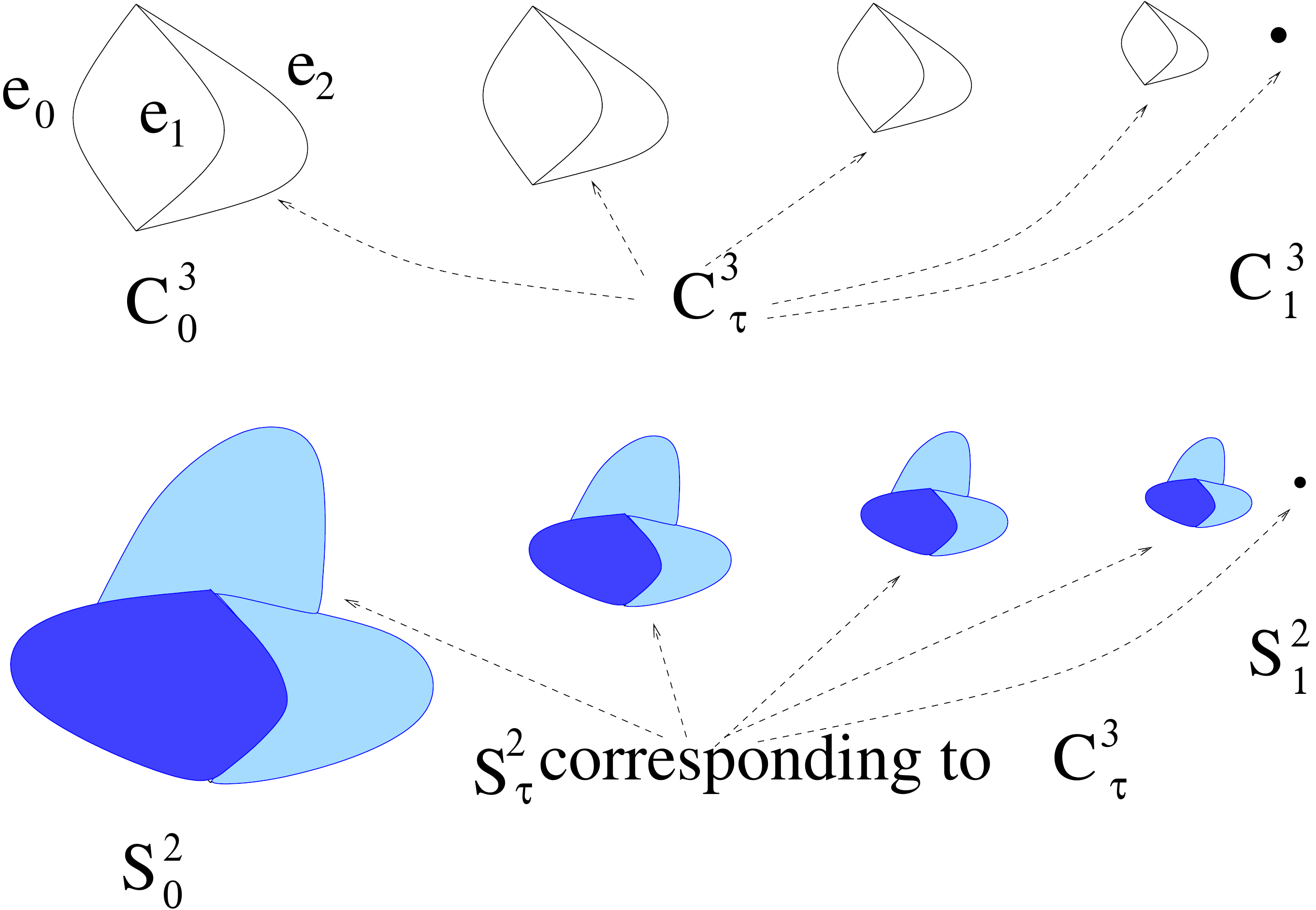}{Deforming a $2$-sphere to a point}{0.5\hsize}

\section{The proof of Theorem ~\ref{TheoremB}.}

In this section we will prove  Theorem ~\ref{TheoremB}. The volume upper bound for the length
of  a ``wide'' geodesic loop will follow from the Filling Radius upper bound for the length 
of this ``wide'' loop combined with the volume upper bound for the Filling Radius.

In [G] Gromov defines the Filling Radius of $M^n$, $FillRad M^n$, as the minimal $r$ such that the image of $M^n$ under the Kuratowsky embedding into 
$L^{\infty}(M^n)$ bounds in its $r$-neighborhood in $L^{\infty}(M^n)$. Here
Kuratowsky embedding is the map that sends each point $x \in M^n$ to the 
distance function $d(x, *)$.  



M. Katz proved that $Fill Rad M^n \leq \frac{d}{3}$, where
$d$ is the diameter of $M^n$, (see [K]), while M. Gromov had 
found the first estimate for the filling radius of
a closed Riemannian manifold in terms of the volume of $M^n$, (see [G]).

\begin{Thm} \label{TheoremC}{\bf [G]}
Let $M^n$ be a closed connected Riemannian manifold. Then 
$FillRadM^n \leq g(n) (vol(M^n))^{\frac{1}{n}}$, where
$g(n)=(n+1)n^n(n+1)!^{\frac{1}{2}}$ and $vol(M^n)$ denotes the
volume of $M^n$.
\end{Thm}

The constant $g(n)$ was improved to $27^n(n+1)!$ by S. Wenger in [W].
It was further improved by A. Nabutovsky in [N] to simply $n$.

In this section we will prove the following

\begin{Thm} \label{TheoremD}
Let $M^n$ be a closed Riemannian manifold.  For each integer 
$a \geq 3$ there exists a non-trivial 
stationary geodesic flower 
of total length $2 (n+1)!^2 a^{(n+1)^3} Fill Rad M^n$ which is a critical 
point of the weighted length functional $L_a$.
\end{Thm}

The last missing ingredient in the proof of Theorem ~\ref{TheoremD} is 
the following  straightforward generalization 
of the Merging Lemma proven in [R4]. It is analogous to the first assertion of Lemma ~\ref{Lemmadiameter}.

Let $K$ be a net that corresponds to the $1$-skeleton of an $m$-simplex. Let
$w_0,...,w_{m+1}$ be its vertices. Consider the set of pairs of vertices 
$S=\{(w_i,w_j) | 0 \leq i < j \leq m+1 \}$ with an alphabetical order. 
That is $(w_{i_1}, w_{j_1}) < (w_{i_2}, w_{j_2})$ if and only if $i_1 < i_2$ or
$i_1=i_2$ and $j_1 < j_2$.  Note there is a one-to-one correspondence between the elements of 
$S$ and the set of edges of $K$. Thus the aplphabetical order on $S$ induces an order on the edges of $K$. 

\begin{Lem} {\bf (Merging Lemma)}
Let $K$ be a net corresponding to the $1$-skeleton of an $m$-simplex 
with the vertices $w_0,....,w_{m+1}$ in a Riemannian manifold $M^n$.  
Let $e_i, i=1,...,\frac{(m+1)(m+2)}{2}$ be the edges,
enumerated in the order that corresponds to the alphabetical order
on the set $S$ above.
Consider the following weighted length functional on $K$:
$L(K)=\Sigma_{j=1}^{\frac{(m+1)(m+2)}{2}}a^{j-1}length(e_j)$, where $a \geq 3$.
Then the only non-trivial critical points of this functional
are minimal geodesic flowers.
\end{Lem}

The proof of the Lemma is presented in [R4] with $a=3$. The proof will 
not be changed if one substitutes an arbitrary integer $a$ that is 
greater than or equal to $3$ in lieu of $3$, thus it will not be presented in this paper. We will, however, present the proof in the case of $m=3$ for 
the sake of readability.

\begin{Pf}{\bf Proof of the Merging Lemma for $m=3$}
 Let us begin by 
demonstrating that a non-degenerate $1$-skeleton of a $3$-simplex cannot
be a critical point of $L_a= \Sigma_{i=1}^6a^{i-1}length (e^i)$, (see fig. ~\ref{flowers2}(a)). 
Let
$v_{ij}$ denote the unit vector tangent to $e_i$ at $w_j$.
Let us consider in turns the four stationarity conditions at each of the 
vertex: $w_0,...,w_3$. 

\noindent (1) The stationarity condition at the vertex $w_0$ impliest that 
$v_{10}+a^2v_{30}+a^3v_{40}=0$;

\noindent (2) The stationarity condition at the vertex $w_1$ implies that
$v_{11}+av_{21}+a^4v_{51}=0$;

\noindent (3) The stationarity condition at the vertex $w_2$ implies that
$av_{22}+a^2v_{32}+a^5v_{62}=0$;

\noindent (4) Finally, the stationarity condition at the last vertex, $w_3$ implies that $a^3v_{43}+a^4v_{53}+a^5v_{63}=0$.

It is obvious that for $a \geq 3$ these conditions cannot be satisfied, unless
each vertex merges with 
some other vertex as, for example is depicted  in the configurations in fig. 
~\ref{flowers2} (b).

Without loss of generality assume that $w_0$ merges with $w_1$. This can only happen
if the length of edge $e_1$ that connects the above vertices decreases to $0$. 
In this case, what can possibly happen with the remaining two vertices $w_2$, and $w_3$?
Since we know that they have to merge with some other vertices, either $w_2$ and $w_3$ merge together, i. e. the length of edge $e_6$ also decreases to $0$, (see fig. ~\ref{flowers2}) or
all of the four vertices join together, and the net becomes a flower.

Thus, it remains to show that the net depicted in fig. ~\ref{flowers2} (b) cannot 
satisfy the stationarity condition. 

In this case the following stationarity condition should be satisfied
at the vertex $w_0=w_1$ 
$av_{02}+a^2v_{03}+a^3v_{12}+a^4v_{13}=0$, (and a similar stationarity condition should be satisfied at $w_2=w_3$. Even if the three edges $e_2,e_3, e_4$ which appear with the lower weights all merge together, it would not be enough to compensate for the weight of $e_5$. Thus the only way the stationarity condition can be satisfied is if all of the vertices merge together.






\end{Pf}

\realfig{flowers2}{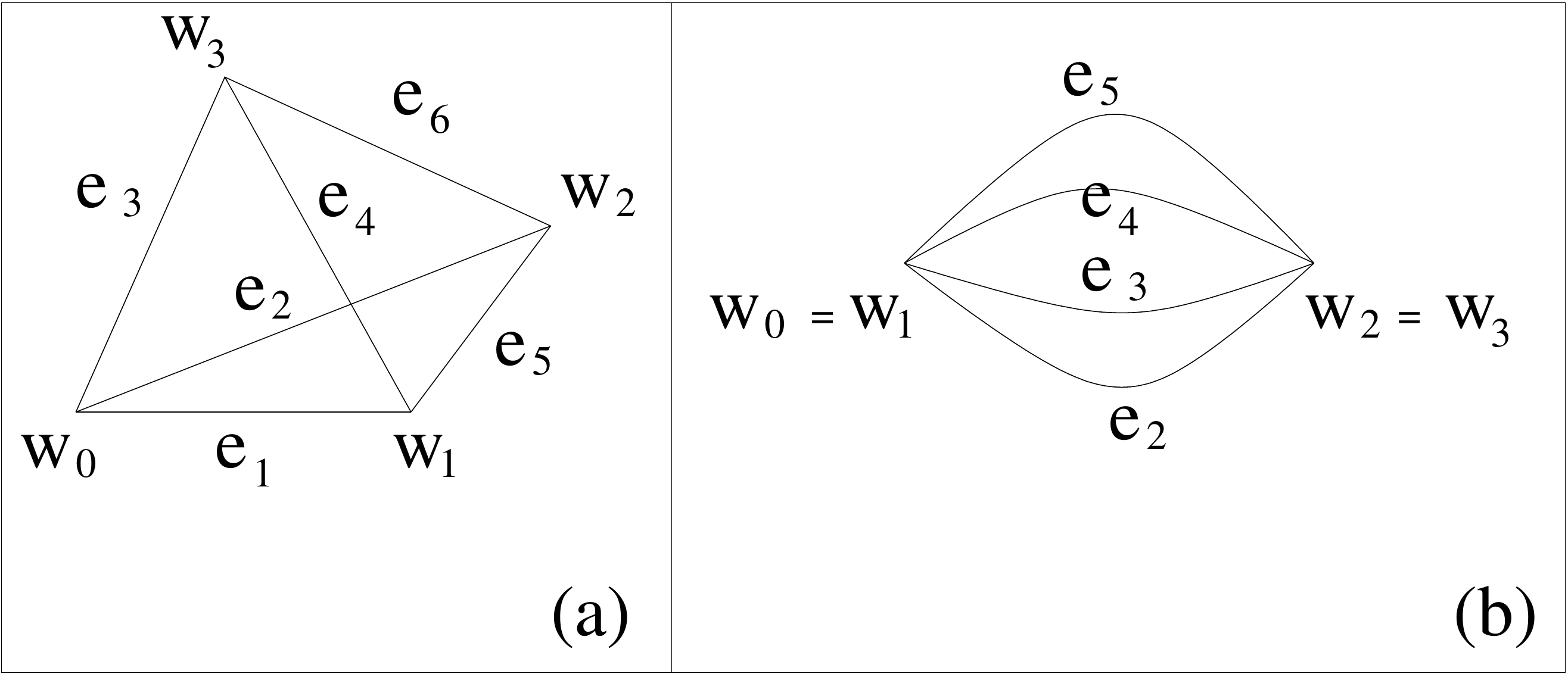}{The graphs below cannot
minimize the weighted length functional $l(\Gamma)=\Sigma_{j=1}^6 a^{j-1}e_j$.}{0.5\hsize}

Theorem ~\ref{TheoremD} combined with Theorem ~\ref{TheoremC}
leads to  the volume bound for the length of the smallest critical point
of $L_a$. By the Merging lemma, this critical point will be a geodesic flower. 
Moreover, it will follow  from the second statement of Lemma ~\ref{Lemmadiameter} that one of the geodesic loops willl have an angle that is within $\epsilon$ from $\pi$. 

We are now ready to present the proof of Theorem ~\ref{TheoremD}.
Let us  outline its  proof. Let us fill $M^n$ by a polyhedron $P^{n+1}$ in 
$L^{\infty}(M^n)$. Let $id: M^n \longrightarrow M^n$ be the identity map on $M^n$. It is impossible to extend the above map to $P^{n+1}$. Thus, attempting to extend the identity map on $M^n$ to $W^{n+1}$ should fail, and 
the required geodesic flower that is a critical point of $L_a$ will be an obstruction to this extension.
The proof is similar to the proof of the analogous statement of [R4].

The difference between the proofs of Theorems ~\ref{TheoremA} and ~\ref{TheoremD} is 
that instead of contracting $m$-cages, we will be contracting 
$1$-skeleta of simplices.  The contraction of $1$-skeleta will generate the 
1-parameter family of (not geodesic) nets. Correspondingly, we can build the spheres and the discs out of these $1$-skeleta. The weighted length functionals used to contract these nets  will be 
$\Sigma_{i=1}^{(n+1)(n+2)/2}a^{i-1} length(e_i)$, where
$e_i$'s are the edges of the $1$-skeleton of an $(n+1)$-dimensional
simplex, satisfying the condition that the edge
$e_1$ is coming out of the vertex $w_0$, edges
$e_2, e_3$ are coming out of the vertex $w_1$, edges
$e_4,e_5,e_6$ are coming out of the vertex $w_3$, etc. Finally, 
the last $(n+1)$ edges are
coming out of the same vertex $w_{n+1}$. As was stated before, the functional 
will force the net to degenerate into a flower, and one of the loops of the flower to become wide.

\begin{Pf}{Proof of Theorem ~\ref{TheoremD}}
Suppose there are no ``small'' critical geodesic flowers satisfying the conditions stipulated in the statement of the theorem.

Let $P^{n+1} \subset L^\infty (M^n)$ be a polyhedron that satisfies: 

\noindent  (1) $\partial P^{n+1} = M^n$, when $M^n$ is orientable, and
$\partial P^{n+1} = M^n$ mod $2$, when $M^n$ is not orientable;

\noindent (2) $\partial P^{n+1} \subset (FillRad M^n + \delta)$-neighborhood 
of $M^n$ for an arbitrarily small $\delta >0$, (see [G]).


Suppose $P^{n+1}$ is  triangulated, and that the diameter
of any simplex in this triangulation is smaller than $\delta$. Note that we obtain a triangulation of $M^n$ by 
restricting the triangulation of $P^{n+1}$.

We will obtain geodesic flowers  of ``small'' length, and with at 
least one ``wide'' geodesic loop as an obstruction to the extension of the identity map $id: M^n \longrightarrow M^n$
to $P^{n+1}$.

This extension  is constructed by induction on the dimension of the skeleta
of $P^{n+1}$.  

As in the proof of Theorem ~\ref{TheoremA}, we will first extend to the {\bf $0$-skeleton}. Assign to  each $\tilde{p}_i
\in P^{n+1}$, the vertex of the triangulation of $P^{n+1}$ the vertex $w_i \in M^n$ which is closest to 
$\tilde{p}_i$. In case there is more than one vertex in the triangulation of $M^n$ that is closest to $\tilde{p}_i$ we can arbitrarily choose any one.   
Therefore, $d(\tilde{p}_i,w_i) \leq FillRad M^n + \delta$. For the sake of the future reference we can call this {\bf step 0} of the extension procedure. 

Next, we will extend to the {\bf $1$-skeleton}, ({\bf step 1}). Let 
$[\tilde{p}_i, \tilde{p}_j] \subset P^{n+1} \backslash M^n$ denote a $1$-simplex
in the triangulation of $P^{n+1}$.
We will assign to it 
a minimal geodesic segment
$[w_i,w_j]$ that connects $w_i$ and $w_j$ of length $\leq 2FillRad M^n +3\delta$. 
Again, in case there are several minimizing geodesics connecting $w_i$, $w_j$, 
we can choose any one of them.

Next we will extend to  the {\bf $2$-skeleton} ({\bf step 2}). Let 
$\tilde{\sigma}^2_{i_0,i_1,i_2}=[\tilde{p}_{i_0},\tilde{p}_{i_1},\tilde{p}_{i_2}]$ be an arbitrary $2$-simplex.  Its boundary is mapped to a closed
curve made out of three geodesic segments that we obtained during step 1 of the extension. 
The total  length of the boundary is at most $\leq 6 FillRad M^n+9\delta$.  If there is a periodic geodesic
of length  at most $6 FillRad M^n+9\delta$, the conclusion of 
Theorem ~\ref{TheoremD} is satisfied. In case there are no periodic geodesics
on $M^n$ satisfying the above bound
we can contract this curve to a point along the curves of smaller length.
Moreover, the absence of ``short'' periodic geodesics implies that
this curve shortening homotopy continuously depends on the initial  curve.  
We will map $\tilde{\sigma}^2_{i_0,i_1,i_2}$ 
to a surface denoted as $\sigma^2_{i_0,i_1,i_2}$, that is generated by the above homotopy.  

Next let us extend to 
the {\bf $3$-skeleton} ({\bf step3}).  
Consider an arbitrary $3$-simplex
$\tilde{\sigma}^3_{i_0,i_1,i_2,i_3}=[\tilde{p}_{i_0},...,\tilde{p}_{i_3}]$.
We know that the boundary of this simplex  is mapped to
the chain: 
$\Sigma_{j=0}^3 (-1)^j \sigma^2_{i_i,...,\hat{i}_j,...,i_3}$ by the previous step of the induction.  Its one skeleton  is  a net that we will  
denote by $K$ defined in Step 1.  Apply the 
 weighted length shortening process for nets. In the absence of the critical geodesic flowers it will be  continuously 
deformed  to a point. The weighted length of $K$ is defined as 
$\Sigma_{j=1}^6 a^{j-1} length (e_j)$, where $e_j$ is an edge of the 
$1$-skeleton.  
(We will not explicitly describe this length shortening process, but as we have mentioned before,
it can be found in [NR1]).

By the Merging Lemma above only geodesic flowers  
can be critical points for 
such a functional.


Next for each time $\tau$ 
one can construct a $2$-dimensional sphere $S^2_\tau$ corresponding to each 
$K_\tau$ obtained by shortening $K=K_0$. The sphere is constructed in a 
way that is analogous to a similar construction 
in the proof of Theorem ~\ref{TheoremA}.  
That is for each $\tau \in [0,1]$ consider $K_\tau$. Consider all triples of the edges
that correspond to the boundary of the face of the three simplex. Apply the curve shortening to each of the closed curves formed by these curves. In the absence of short periodic geodesics the process will converge to a point. Thus, we will obtain the $2$-discs that these homotopies will generate. Glue those discs  together along the common edges to obtain $S^2_\tau$. 
This $1$-parameter family
of $2$-spheres can be regarded as a $3$-disc that we will denote
as $\sigma^3_{i_0,...,i_3}$.  We will assign it to 
$\tilde{\sigma}^3_{i_0,...,i_3}$.  

If we  continue this inductive procedure 
until we reach the $(n+1)$-skeleton of $P$. We will obtain  
a singular chain on $M^n$, that has the fundamental class $[M^n]$ as its
boundary, and therefore, arrive at a contradiction.


 

  Suppose we have extended the identity map
$id: M^n \longrightarrow M^n$ to the $k$-skeleton of $P^{n+1}$
and now want to extend it to the $(k+1)$-skeleton.  Take
an arbitrary $(k+1)$-simplex of $P^{n+1}$. Let us denote by $N$ the image of its
$1$-skeleton. To extend the map to 
this simplex, we will construct a $(k+1)$-dimensional disc that
fills $N$.  $N$ consists of $\frac{(k+2)(k+1)}{2}$ edges.   
Each edge has a weight assigned to it
as in the Merging Lemma above.  Assuming there is no ``small'' geodesic
flowers, $N$ can be deformed to a point along the $1$-parameter
family of nets $N_\tau, \tau \in [0,1]$, where the weighted length
of $N_\tau$ decreases with $\tau$.  Next for each $\tau$ we
construct a sphere $S^k_\tau$ that fills $N_\tau$.  It is consructed
by constructing $k$-dimensional discs and gluing them as in the
boundary of a $(k+1)$-dimensional simplex.  In order to construct
those discs, we consider subnets that are obtained by ignoring
a vertex and all the edges that are coming out of this vertex.

The maximal length of each edge in a subnet is $\leq 2 Fill Rad M^n \frac{(k+2)(k+1)}{2}a^{\frac{(k+2)(k+1)-2}{2}}$,
i.e. the maximal number of edges in $N_\tau$ times the maximal 
weighted length of the edges in $N_\tau$.  Thus, the total 
weighted length of the whole subnet is 
$\leq 2 Fill Rad M^n \frac{(k+1)k(k+2)(k+1)}{2^2} a^{\frac{(k+1)k-2}{2}}
a^\frac{(k+2)(k+1)-2}{2}$.  Proceeding in the manner starting
from the $(n+1)$-skeleton of $P^{n+1}$ we would obtain a bound
of $2 Fill Rad M^n \frac{(n+2)!(n+1)!}{2^n}
a^{\Sigma_{k=0}^n\frac{(k+2)(k+1)-2}{2}} \leq 2 Fill Rad M^n (n+1)!^2
a^{(n+1)^3}$.

Note also that the maximal number of geodesic loops in the geodesic
flower can be estimated by the number of edges in $N$, but taken with
multiplicities that correspond to weights, thus it is bounded by
the sum $\Sigma_{j=1}^{\frac{(n+2)(n+1)}{2}}a^{j-1} \leq a^{(n+1)^2}$.
\end{Pf}

\begin{Pf}{Proof of Theorem ~\ref{TheoremB}} We combine the asertions of 
Theorem ~\ref{TheoremD} and Lemma 1.1. The desired assertion follows when we substitute $a=\max \{ \lceil {1 \over \sin{ \epsilon \over 2}} \rceil+1, 3 \}+ \delta$ in the inequality of Theorem ~\ref{TheoremD}, and then take the limit as $\delta \longrightarrow 0$.
\end{Pf}

\noindent{ \bf Acknowledgments.} The author gratefully acknowledges the 
partial support by Natural Sciences and Engineering Research Council
(NSERC) University Faculty Award her work on the present paper.
This paper was partially written during the author's visit of
the Institute for Advanced Studies, Princeton.  
She would like to thank
the Institute for Advanced Studies for its kind hospitality.


\normalsize

\bigskip
\begin{tabbing}
\hspace*{7.5cm}\=\kill
R. ~Rotman\\
Department of Mathematics                  \\
University of Toronto                     \\
Toronto, Ontario M5S 2E4                  \\      
Canada                                    \\ 
e-mail: rina@math.toronto.edu             \\

\end{tabbing}

\end{document}